\pdfoutput=1

\documentclass{amsart}
\usepackage{tikz}

\usepackage{amssymb}
\usepackage{amsmath}
\usepackage{amsfonts}
\usepackage{url}
\usepackage{enumerate}
\usepackage{graphicx}
\usepackage{mathtools}
\usepackage{mathrsfs}
\usepackage{newfloat}
\DeclareFloatingEnvironment[
  name = Procedure
]{procedure}
\numberwithin{procedure}{section}

\usepackage{hyperref}
\hypersetup{
    colorlinks  =   true,
    linkcolor   =   blue,
    filecolor   =   blue,
    urlcolor    =   blue,
    citecolor   =   blue,
}

\usepackage{algpseudocode}

\newcommand{\erase}[1]{}

\theoremstyle{plain}
\newtheorem{theorem}{Theorem}[section]       
\newtheorem{corollary}[theorem]{Corollary}   

\theoremstyle{definition}
\newtheorem{definition}[theorem]{Definition} 

\theoremstyle{remark}
\newtheorem{remark}[theorem]{Remark}         

\numberwithin{equation}{section}
\numberwithin{table}{section}
\numberwithin{figure}{section}

\newcommand{\qedbox}{\hfill {$\Box$}}

\newcommand{\CC}{\mathord{\mathbb{C}}}

\newcommand{\FF}{\mathord{\mathbb{F}}}

\newcommand{\PP}{\mathord{\mathbb{P}}}

\newcommand{\ZZ}{\mathord{\mathbb{Z}}}

\newcommand{\AAA}{\mathord{\mathcal{A}}}

\newcommand{\CCC}{\mathord{\mathcal{C}}}

\newcommand{\LLL}{\mathord{\mathcal{L}}}

\newcommand{\OOO}{\mathord{\mathcal{O}}}

\newcommand{\UUU}{\mathord{\mathcal{U}}}

\newcommand{\XXX}{\mathord{\mathcal{X}}}

\newcommand{\SSSS}{\mathord{\mathfrak{S}}}    

\newcommand{\inj}{\hookrightarrow}
\providecommand{\twoheadrightarrow}{\mathrel{\rightarrow\!\!\!\!\rightarrow}}
\newcommand{\surj}{\twoheadrightarrow}

\newcommand{\isom}{\xrightarrow{\raise -3pt \hbox{\scriptsize $\sim$} }}

\newcommand{\set}[2]{\{\,{#1}\mid {#2} \,\}}           
\newcommand{\bigset}[2]{\left\{\; {#1} \; \left\vert \; {#2} \;  \right.\right \}}

\newcommand{\angs}[1]{\langle {#1}  \rangle}           

\newcommand{\sprime}{^{\prime}}

\newcommand{\spar}[1]{^{(#1)}}
\newcommand{\spprime}{^{\prime\prime}}

\newcommand{\sperp}{^{\perp}}

\newcommand{\semidirectproduct}{\rtimes}    

\DeclareMathOperator{\Coker}{Coker}

\newcommand{\OG}{\mathord{\mathrm{O}}}

\newcommand{\pione}{\pi_1}  

\newcommand{\intf}[1]{\langle #1 \rangle}

\newcommand{\intfnull}{\intf{\phantom{a}, \phantom{a}}}
\newcommand{\diag}{\mathord{\mathrm{diag}}}

\newcommand{\Grass}{\mathord{\mathrm{Grass}}}
\newcommand{\Pic}{\mathord{\mathrm{Pic}}}

\newcommand{\nbsim}{\mathord{\sim}}
\newcommand{\nbsimd}{\nbsim_d} 
\newcommand{\nbsimc}{\nbsim_c}
\newcommand{\nbsimt}{\nbsim_t}
\newcommand{\bsimd}{\mathbin{\nbsimd}}
\newcommand{\bsimc}{\mathbin{\nbsimc}}
\newcommand{\bsimt}{\mathbin{\nbsimt}}

\begin{document}

\title[Zariski pairs on cubic surfaces]
{Zariski pairs  on   cubic surfaces}

\author[Ichiro  Shimada]{Ichiro Shimada}
\address{Mathematics Program, Graduate School of Advanced Science and Engineering, Hiroshima University, 
1-3-1 Kagamiyama, Higashi-Hiroshima, 739-8526 JAPAN}
\email{ichiro-shimada@hiroshima-u.ac.jp}

\begin{abstract}
A line arrangement on a smooth cubic surface is 
a subset of the set of  lines lying on the surface.
We define the notion of 
Zariski pairs of line arrangements on general cubic surfaces,
and provide a complete classification of these Zariski pairs.
\end{abstract}


\thanks{This work was supported by JSPS KAKENHI,  
~20K20879, ~20H01798 and~23H00081.}
%
%
%
%
%
%

\maketitle
\section{Introduction}\label{sec:introduction}
We work over the complex number field $\CC$.
Cayley and Salmon showed in 1849 that every smooth cubic surface contains exactly $27$ lines.
The configuration of these $27$ lines  is a beautiful historical topic of algebraic geometry.
In this paper, we investigate this configuration from the viewpoint of
\emph{Zariski pairs}.
\par
By a \emph{plane curve},
we mean  a reduced, possibly reducible,  projective plane curve.
We say that a pair $(C_1, C_2)$ of plane curves is a \emph{Zariski pair} 
if $C_1$ and $C_2$ have the same combinatorial type of singularities, 
but have different embedding topologies in the projective plane.
This notion of Zariski pairs  was formulated in Artal~Bartolo's seminal paper~\cite{Artal1994},
in which he investigated a pair of $6$-cuspidal sextics 
discovered by Zariski in 1929,   
and presented some new examples of Zariski pairs.
Since then,
many authors have studied Zariski pairs of plane curves from various points of view.
See, for example, the survey~\cite{Survey2008}.
\par
We introduce the notion of  \emph{Zariski pairs of line arrangements
on  general cubic surfaces}.
\begin{definition}
A point $Q$ of a smooth cubic surface $X$ is called an \emph{Eckardt point}
if three lines on $X$ pass through $Q$.
\end{definition}
A general cubic surface has no Eckardt points.
Let $X\subset \PP^3$ be a smooth cubic surface with no Eckardt points, 
and let $L(X)$ denote the set of lines on $X$.
We describe the configuration of lines on $X$ by the intersection form 
\[
\intf{\ell, \ell\sprime}:=
\begin{cases}
-1 &\textrm{if $\ell=\ell\sprime$,}\\
0 &\textrm{if $\ell\ne \ell\sprime$,  and $\ell$ and $\ell\sprime$ are disjoint,}\\
1 &\textrm{if $\ell\ne \ell\sprime$,   and $\ell$ and $\ell\sprime$ intersect, }\\
\end{cases}
\]
for $\ell, \ell\sprime\in L(X)$.
\begin{definition}
A \emph{line arrangement on a general cubic surface} is a pair $[S, X]$ 
of a smooth cubic surface $X$ with no  Eckardt points and a subset $S$ of   $L(X)$.
In this situation, 
we say that $S$ is a \emph{line arrangement on $X$}.
We denote by $\AAA$ the set of line arrangements on  general cubic surfaces.
\end{definition}
We introduce three equivalence relations $\nbsimd$, $\nbsimc$, and $\nbsimt$ on $\AAA$.
\begin{definition}
Let $[S, X]$ and $[S\sprime, X\sprime]$ be elements of $\AAA$.
\begin{itemize}
\item 
We say that $[S, X]$ and $[S\sprime, X\sprime]$ are \emph{deformation equivalent}
and write 
\[
[S, X]\bsimd[S\sprime, X\sprime]
\]
if there exists a continuous family 
$\XXX:=\set{X_t}{t\in [0,1]}$ of  smooth cubic surfaces with no Eckardt points
connecting $X=X_0$ and $X\sprime=X_1$
such that 
$S$ is deformed continuously to $S\sprime$ 
along $\XXX$.
We denote by $[S, X]_d$ the equivalence class containing $[S, X]$ under the equivalence relation $\nbsimd$.
\item 
We say that $[S, X]$ and $[S\sprime, X\sprime]$ have the \emph{same embedding topology}  
and write 
\[
[S, X]\bsimt [S\sprime, X\sprime]
\]
if there exists a homeomorphism $X\isom X\sprime$
that maps the union $\Lambda(S)\subset X$ of lines in $S$ 
to the union $\Lambda(S\sprime)\subset X\sprime$ of lines in $S\sprime$.
We denote by $[S, X]_t$ the equivalence class containing $[S, X]$ under the equivalence relation $\nbsimt$.
\item
We say that  $[S, X]$ and $[S\sprime, X\sprime]$ have the  \emph{same combinatorial type}
and write 
\[
[S, X]\bsimc [S\sprime, X\sprime]
\]
if there exists a bijection  between $S$ and $ S\sprime$ 
that preserves the intersection form $\intfnull$.
We denote by $[S, X]_c$ the equivalence class containing $[S, X]$ under the equivalence relation $\nbsimc$.
\end{itemize}
\end{definition}
It is obvious that we have the following implications:
\begin{equation*}\label{eq:implications}
[S, X]\bsimd[S\sprime, X\sprime]\;\;\Longrightarrow\;\;
[S, X]\bsimt[S\sprime, X\sprime]\;\;\Longrightarrow\;\;
[S, X]\bsimc[S\sprime, X\sprime].
\end{equation*}
Therefore we have natural surjections
\[
 \AAA/\nbsimd \;\;\surj\;\;  \AAA/\nbsimt \;\;\surj\;\;  \AAA/\nbsimc.
\]
Following the definition of Zariski pairs of plane curves,
we make the following: 
\begin{definition}
We  say that two equivalence classes $[S, X]_d$ and $[S\sprime, X\sprime]_d$  
form 
a \emph{Zariski pair  of line arrangements on general cubic surfaces} (a Zariski pair in $\AAA$,  for short)
if $[S, X]$ and $[S\sprime, X\sprime]$ have the same combinatorial type, 
but have different embedding topologies.
\end{definition}
We choose and fix a smooth cubic surface $X$ with no Eckardt points, 
and denote by $\AAA_X:=2^{L(X)}$ the set of line arrangements on $X$.
Note that smooth cubic surfaces with no Eckardt points are parameterized by
a Zariski open subset $\UUU^0$ of the projective space $\PP_*(H^0(\PP^3, \OOO(3)))$
parameterizing all cubic surfaces.
Since $\UUU^0$ is 
 connected, 
 the inclusion $\AAA_X\inj \AAA$
induces a bijection 
\[
 \AAA_X/\nbsimd \;\;\cong\;\;  \AAA/\nbsimd.
\]
Since $\AAA_X$ is  finite, 
we can regard Zariski pairs in $\AAA$
as a \emph{toy model} of
classical Zariski pairs of plane curves.
In fact, 
 we can enumerate all Zariski pairs in $\AAA$ 
by a brute force method.
This complete list  is the main result of this note.
 \par
To distinguish  embedding topologies, 
we use the lattice structure 
on  the middle cohomology group $H^2(X, \ZZ)$ of the smooth cubic surface $X$.
The cup-product $\intfnull$ makes 
 $H^2(X, \ZZ)$  a unimodular lattice of rank $7$.
 For a line arrangement $S$ on $X$,
 let $H(S)\subset H^2(X, \ZZ)$ denote the submodule generated by 
 the classes of lines in $S$,
 and we put 
 \[
 H(S)\sperp:=\set{x\in H^2(X, \ZZ)}{\intf{x, y}=0\;\;\textrm{for all}\;\; y\in H(S)}.
 \]
Remark that $H(S)$ and $H(S)\sperp$ are topological invariants of the pair $[S, X]$.
 \par
  Recall that a lattice $M$ is said to be \emph{even} if $\intf{x, x}\in 2\ZZ$ holds for all $x\in M$,
  and to be  \emph{odd} otherwise.
Our result is as follows:
\begin{theorem}\label{thm:main}
There exist exactly two Zariski pairs 
\[
([S_1, X]_d,\;\;  [S_2, X]_d)
\quad{\rm and}\quad
([T_1, X]_d ,\;\;  [T_2, X]_d)
\]
of line arrangements  on general cubic surfaces.
\begin{enumerate}[{\rm(1)}]
\item The combinatorial type of $S_i$ is as follows.
We have $|S_i|=5$, and any distinct lines $\ell, \ell\sprime\in S_i$ are disjoint.
The embedding topologies of $S_1$ and $S_2$ are distinguished by  the fact that 
 $H(S_1)\sperp$ is odd,
whereas $H(S_2)\sperp$ is even.
\item The combinatorial type of $T_i$ is as follows.
We have $|T_i|=6$, and, 
 for $\mu\ne \nu$, 
\[
\intf{\ell_\mu, \ell_\nu}=\begin{cases}
1 & \textrm{if $\mu=0$ or $\nu=0$,} \\
0 & \textrm{if $\mu\ne 0$ and $\nu\ne 0$} \\
\end{cases}
\]
holds under a suitable numbering $\ell_0, \dots, \ell_5$ of the elements of $T_i$.
The embedding topologies of $T_1$ and $T_2$ are distinguished by the following:
\[
H_1(X\setminus \Lambda(T_1), \ZZ)\cong \ZZ/2\ZZ, \quad H_1(X\setminus \Lambda(T_2), \ZZ)=0.
\]
\end{enumerate}
\end{theorem}
The main ingredient of 
the proof  is  the result of Harris~\cite{Harris1979}
on the Galois group  of the $27$ lines on a  cubic surface $X$.
We write the Galois action on $L(X)$ explicitly, and 
calculate  the orbit-decomposition 
of   $\AAA_X=2^{L(X)}$.
Comparing the combinatorial types and the embedding topologies of these orbits,
we obtain Theorem~\ref{thm:main}.
\par
A cubic surface is a del Pezzo surface of degree $3$, and 
the Galois group of its $27$ lines  is isomorphic to the Weyl group $W(E_6)$ of type $E_6$.
In~\cite{Shimada2022}~and~\cite{Shimada2025},
we investigated Zariski multiples associated 
with del Pezzo surfaces of degree~$2$ and~$1$, 
using the Galois actions of $W(E_7)$ and $W(E_8)$, respectively.
For general methods of distinguishing embedding topologies via  lattices, see~\cite{Shimada2010}.
\par
For the actual computation, we used {\tt GAP}~\cite{GAP4}.
In~\cite{WE6compdata}, we present a detailed computation data. 
\par
\medskip
{\bf Convention.}
The orthogonal group $\OG(M)$ of a lattice $M$ acts on $M$ from the right.
The symmetric group $ \SSSS(T)$ of a finite set $T$ also acts on $T$ from the right.
\section{The $27$ lines on a cubic surface}\label{sec:27}
In this section, 
we recall some basic facts about cubic surfaces and review the result of Harris~\cite{Harris1979},
which is  reproduced in \cite[Section 3.3]{Shimada2025} by a simpler method.
For a general theory  of  cubic surfaces, we refer the reader to Demazure~\cite{Demazure1980} 
or to   Dolgachev~\cite[Chapter 9]{Dolgachev2012}.
\subsection{Action of $W(E_6)$ on the $27$ lines}
Let $P_1, \dots, P_6$ be general six points of $\PP^2$,
and let $X\to \PP^2$
be the blowing-up at $P_1, \dots, P_6$.
For a divisor $D$ on $X$, let $[D]\in H^2(X, \ZZ)$ denote its class.
Then $D\mapsto [D]$ induces an isomorphism
from the Picard group $\Pic \, X$  with the intersection pairing to 
$H^2(X, \ZZ)$ with the cup-product $\intfnull$.
From now on, we identify $\Pic\,(X)$ with $ H^2(X, \ZZ)$.
Let $h\in H^2(X, \ZZ)$ be the class of the pull-back of a line on $\PP^2$,
and let $e_i:=[E_i]$ be the class of the exceptional curve $E_i$  over $P_i$ for $i=1, \dots, 6$.
The lattice $ H^2(X, \ZZ)$ is of rank $7$ with a basis $h, e_1, \dots, e_6$,
with respect to which the Gram matrix is given by the diagonal matrix 
\begin{equation*}\label{eq:diag}
 \diag\, (1, -1,-1,-1,-1,-1,-1).
 \end{equation*}
We express elements of $H^2(X, \ZZ)$ as vectors in terms of this basis.
The class of the anti-canonical divisor $-K_X$ is 
\[
[-K_X]=(3,-1,-1,-1,-1,-1,-1).
\]
We have $\intf{-K_X, -K_X}=3$,
and the complete linear system $|-K_X|$  embeds
$X$ into $\PP^3$ as a smooth cubic surface.
We denote by $K$ the sublattice of $H^2(X, \ZZ)$ generated by $[-K_X]$,
and by $V$ the orthogonal complement of $K$ in $H^2(X, \ZZ)$.
Then $V$ is a negative-definite root lattice of type $E_6$.
Indeed, the $(-2)$-vectors 
\begin{eqnarray*}
r_{ 1 } &:=& ( -1, 0, 0, 0, 1, 1, 1 ) , \\ 
r_{ 2 } &:=& ( 0, 1, -1, 0, 0, 0, 0 ) , \\ 
r_{ 3 } &:=& ( 0, 0, 1, -1, 0, 0, 0 ) , \\ 
r_{ 4 } &:=& ( 0, 0, 0, 1, -1, 0, 0 ) , \\ 
r_{ 5 } &:=& ( 0, 0, 0, 0, 1, -1, 0 ) , \\ 
r_{ 6 } &:=& ( 0, 0, 0, 0, 0, 1, -1 ) ,
\end{eqnarray*}
constitute  a basis of the lattice  $V$,  and  form the dual graph
\begin{equation}\label{eq:E6}
\raise -.7cm \hbox{
\begin{tikzpicture}[ x=12mm, y=12mm]
\coordinate (e1) at (0,.7);
\coordinate (e2) at (-2,0);
\coordinate (e3) at (-1,0);
\coordinate (e4) at (0,0);
\coordinate (e5) at (1,0);
\coordinate (e6) at (2,0);
\draw (e1)--(e4);
\draw (e2)--(e3);
\draw (e3)--(e4);
\draw (e4)--(e5);
\draw (e5)--(e6);
\fill (e1) circle [radius=.1];
\fill (e2) circle [radius=.1];
\fill (e3) circle [radius=.1];
\fill (e4) circle [radius=.1];
\fill (e5) circle [radius=.1];
\fill (e6) circle [radius=.1];
\node at (e1) [right]{\;$r_1$};
\node at (-2,-.1)[below]{$r_2$};
\node at (-1,-.1)[below]{$r_3$};
\node at (-0,-.1)[below]{$r_4$};
\node at (1,-.1)[below]{$r_5$};
\node at (2,-.1)[below]{$r_6$};
\end{tikzpicture}
},
\end{equation}
which is the Dynkin diagram of type $E_6$.
\begin{remark}
This collection of 
 $(-2)$-vectors is chosen solely  for the group-theoretic computation of $W(E_6)$ below,
and carries no geometric significance.
\end{remark}
Hence we have 
\[
\OG(V)=W(E_6)\semidirectproduct\angs{g_0},
\]
where $W(E_6)\subset \OG(V)$ is the Weyl group of type $E_6$ generated by the reflections 
\[
\sigma_{\nu}\colon x\mapsto x+\intf{x, r_{\nu}} r_{\nu}\qquad(\nu=1, \dots, 6)
\]
with respect to the roots 
$r_1, \dots, r_6\in V$, and $g_0$ is the involution of $V$ given by
\[
r_1\leftrightarrow r_1,
\quad
r_2\leftrightarrow r_6,
\quad
r_3\leftrightarrow r_5,
\quad
r_4\leftrightarrow r_4,
\]
which corresponds to the automorphism of the graph~\eqref{eq:E6}.
%
\par
By~\cite[Proposition 3.1]{Shimada2025}, 
we have 
\begin{equation*}\label{eq:WE6isom1}
W(E_6)=\bigset{g\in \OG(V)}{\parbox{7cm}{the isometry $g$ extends to an isometry $\tilde{g}$ of $H^2(X, \ZZ)$ that acts on  $K$ trivially}}, 
\end{equation*}
and hence 
the mapping  $g\mapsto \tilde{g}$ gives rise to an isomorphism 
\begin{equation}\label{eq:WE6isom2}
W(E_6)\isom \set{\tilde{g}\in \OG(H^2(X, \ZZ))}{[-K_X]^{\tilde{g}}=[-K_X]}.
\end{equation}
Since  the set $L(X)$ of lines on $X$ is embedded into $H^2(X, \ZZ)$ by $\ell\mapsto [\ell]$, 
we can compute the permutation action 
\begin{equation}\label{eq:WE6SSSSL}
W(E_6) \to \SSSS(L(X))
\end{equation}
explicitly as follows.
\par
For $i=1, \dots, 6$,
let $\ell[i]$ denote the exceptional curve $E_i$ over $P_i$.
For $i, j$ with $1\le i<j\le 6$, 
let $\ell[ij]$ denote the strict transform of the line on $\PP^2$ passing through $P_i$ and $P_j$.
For $k=1, \dots, 6$,
let $\ell[\bar{k}]$ denote the strict transform of the conic on $\PP^2$ 
passing through the $5$ points in $\{P_1, \dots, P_6\}\setminus \{P_k\}$.
The set $L(X)$ consists of these smooth rational curves.
Their classes are
\[
[\ell[i] ]= e_i, \quad
[\ell[ij] ]= h-e_i-e_j, \quad
[\ell[\bar k]]=2h-(e_1+\dots+e_6)+e_k.
\]
%
We number the elements of $L(X)=\{\ell_1, \dots, \ell_{27}\}$ as follows:
\begin{equation}\label{eq:numbering}
\renewcommand{\arraystretch}{1.25}
\begin{array}{lll}
&&
\ell_{1}:=\ell[1], \;  \dots, \; \ell_{6}:=\ell[6],\\
&&
\ell_{7}:=\ell[12],\;
\ell_{8}:=\ell[13],\;
\ell_{9}:=\ell[14],\;
\ell_{10}:=\ell[15],\;
\ell_{11}:=\ell[16],\\
&&\ell_{12}:=\ell[23],\;
\ell_{13}:=\ell[24],\;
\ell_{14}:=\ell[25],\;
\ell_{15}:=\ell[26], \\
&&
\ell_{16}:=\ell[34],\; 
\ell_{17}:=\ell[35],\; 
\ell_{18}:=\ell[36],\; \\
&&
\ell_{19}:=\ell[45],\; 
\ell_{20}:=\ell[46],\; 
\ell_{21}:=\ell[56], \\
&&
\ell_{22}:=\ell[\bar 1], \;  \dots, \;  \ell_{27}:=\ell[\bar 6].
\end{array}
\end{equation}
We let $\tau\in \SSSS_{27}$  act on $L(X)$ as 
$(\ell_{i})^{\tau}:=\ell_{(i^{\tau})}$.
Then the reflections $\sigma_{\nu}\in W(E_6)$  act on $L(X)$ by the following permutations:
%
%
\begin{equation}\label{eq:ssperms}
\renewcommand{\arraystretch}{1.25}
\begin{array}{ccl}
\sigma_1 &\mapsto& ( 4,21)( 5,20)( 6,19)( 7,24)( 8,23)(12,22), \\
\sigma_2 &\mapsto& ( 1, 2)( 8,12)( 9,13)(10,14)(11,15)(22,23), \\
\sigma_3 &\mapsto& ( 2, 3)( 7, 8)(13,16)(14,17)(15,18)(23,24), \\
\sigma_4 &\mapsto& ( 3, 4)( 8, 9)(12,13)(17,19)(18,20)(24,25), \\
\sigma_5 &\mapsto& ( 4, 5)( 9,10)(13,14)(16,17)(20,21)(25,26), \\
\sigma_6 &\mapsto& ( 5, 6)(10,11)(14,15)(17,18)(19,20)(26,27). 
\end{array}
\end{equation}
\subsection{Monodromy action on the $27$ lines}
All cubic surfaces are parameterized by the projective space $\PP^{19}=\PP_*(H^0(\PP^3, \OOO(3)))$.
For $t\in \PP^{19}$, let $X_t\subset \PP^3$ denote the corresponding cubic surface.
We put
\[
\UUU:=\set{t\in \PP^{19}}{\textrm{$X_t$ is smooth}},
\quad
\UUU^0:=\set{t\in\UUU}{\textrm{$X_t$ has no Eckardt points}},
\]
which are Zariski open subsets of $\PP^{19}$.
We then put
\[
\LLL:=\set{(t, \ell)}{\ell \subset X_t}\;\;\subset\;\; \UUU\times \Grass(\PP^1, \PP^3),
\]
where $ \Grass(\PP^1, \PP^3)$ is the Grassmannian variety of lines in $\PP^3$.
The first projection $\pi_{\LLL}\colon \LLL\to \UUU$
is an \'etale covering of degree $27$,
and the fiber $L_t$ of $\pi_{\LLL}$ over $t\in \UUU$ is the set $L(X_t)$ of lines on the cubic surface $X_t$.
Let $b\in \UUU\sp0$ be the point such that $X_b$ is the cubic surface $X$ fixed in the previous subsection.
We have $X_b=X$ and $L_b=L(X)$.
Harris~\cite{Harris1979} proved the following.
See also~\cite[Section 3.3]{Shimada2025}
for a simpler proof.
\begin{theorem}[Harris~\cite{Harris1979}]\label{thm:Harris}
The image of the monodromy action 
\[
\mu_L \colon \pione (\UUU, b) \longrightarrow \SSSS(L_b)
\]
associated with  $\pi_{\LLL}$
is equal to the image of the homomorphism~\eqref{eq:WE6SSSSL}.
\qedbox
\end{theorem}
The inclusion $\UUU^0\inj \UUU$
induces a surjective homomorphism
$\pione(\UUU^0, b)\surj \pione(\UUU, b)$.
Therefore, for $S_1, S_2\in \AAA_{X_b}=2^{L(X_b)}$,
 we see that $ [S_1, X_b]\bsimd [S_2, X_b]$ holds 
 if and only if 
 $S_1$ and $S_2$ belong to the same $W(E_6)$-orbit
under the action of $W(E_6)$ on $2^{L_b}$ induced by~\eqref{eq:ssperms}
via the identification  $L_b\cong \{1, \dots, 27\}$ given by~\eqref{eq:numbering}.
  \begin{remark}
  In fact, Harris~\cite{Harris1979} proved that the image of $\mu_L$ is isomorphic to $\OG^-(6, \FF_2)$,
  which is isomorphic to $W(E_6)$.
   \end{remark}
   %
%
\section{Orbit decomposition and Zariski pairs}\label{sec:orbits}
\subsection{$W(E_6)$-orbits}
Recall that $X=X_b$.
We calculate the orbit decomposition of $\AAA_X=2^{L(X)}$
under the action of  $W(E_6)$. 
By the numbering~\eqref{eq:numbering},
a line arrangement on $X$ is 
expressed as a subset of $\{1,\dots, 27\}$.
We  write a line arrangement $S\subset L(X)$
as an \emph{increasing} sequence $[s_1, \dots, s_n]$ of integers in  $\{1,\dots, 27\}$.
In particular,
for $S=[s_1, \dots, s_n]$ and $\gamma \in \SSSS_{27}$,
we denote by $S^{\gamma}$ the increasing sequence of integers
obtained by \emph{sorting} the set $\{s_1^{\gamma}, \dots, s_n^{\gamma}\}$.
Let 
$\CCC_n\subset \AAA_X=2^{L(X)}$ be the set of line arrangements consisting of $n$ lines.
We introduce the lexicographic order $\prec$ on each $\CCC_n$;
that is, for distinct elements $S\spar{0}=[s_1\spar{0}, \dots, s_n\spar{0}]$ 
and $S\spar{1}=[s_1\spar{1}, \dots, s_n\spar{1}]$ of $\CCC_n$,
we have $S\spar{0}\prec S\spar{1}$ if and only if  $s_i\spar{0}< s_i\spar{1}$ for  the smallest index $i$ such that $s_i\spar{0}\ne s_i\spar{1}$.
%
A line arrangement $S\subset L(X)$ is said to be \emph{minimal} 
if $S$ is minimal with respect to $\prec$
in the orbit 
\[
o(S):=\set{S^\gamma}{\gamma\in W(E_6)}.
\]
%
Every  $W(E_6)$-orbit  in $\AAA_X$ contains a unique minimal element.
Note that, if a sequence $S:=[s_1, \dots, s_n]$ is minimal,
then any initial subsequence $[s_1, \dots, s_m]$ of $S$ is also minimal.
Using this property and employing {\tt GAP}~\cite{GAP4}, 
we obtain the complete list of minimal representatives of the $W(E_6)$-orbits.
The result is presented  in Table~\ref{table:numborbits}.
\begin{table}
{\small
\begin{eqnarray*}
&&\begin{array}{c|cccccccccccccc}
n &0&1&2&3&4 & 5 & 6 & 7 & 8 & 9 &10 &11 &12 & 13\\
\hline
|\textrm{orbits}| & 1& 1& 2& 4& 8 &18 & 39 &73 &135 & 234 &363&509 & 641 & 715
\end{array}
\\
&&
\begin{array}{c|cccccccccccccc}
n &14&15 & 16 & 17 & 18 & 19 & 20 &21 & 22& 23 &24 &25& 26 &27\\
\hline
|\textrm{orbits}| & 715 & 641 & 509 & 363 &234 & 135 & 73 &39& 18& 8& 4& 2& 1& 1
\end{array}
\end{eqnarray*}
}
\caption{Number  of orbits $|\CCC_n/\nbsimd|$}\label{table:numborbits}
\end{table}
In particular, we have 
\[
| \AAA_X/\nbsimd |\;=\;\ |\AAA/\nbsimd|\;=\; 5486.
\]
\subsection{Zariski pairs}
We compare the combinatorial types of  all  $W(E_6)$-orbits in $\AAA_X$. 
It turns out that 
the natural surjection $\AAA_X/\nbsimd \surj \AAA_X/\nbsimc$
has exactly two fibers of size $>1$.
Each of these fibers  is of size $2$, 
and the two elements in each of these fibers have different embedding topologies,
as is shown in Sections~\ref{subsubsec:5}~and~\ref{subsubsec:6} below.
As a corollary, we obtain the following:
\begin{corollary}
The equivalence relations~$\nbsimd$ and~$\nbsimt$
on $\AAA$ are the same.
\qedbox
\end{corollary}
\subsubsection{A pair of arrangements of $5$ lines}\label{subsubsec:5}
Let $o(S_1)$ and $o(S_2)$ be the $W(E_6)$-orbits whose minimal representatives are
\[
S_1=[1,2,3,4,5]\quad{\rm and}\quad S_2=[1,2,3,4,21], 
\]
respectively.
We have $|o(S_1)|=432$ and $|o(S_2)|=216$.
Each of these  arrangements consists of disjoint $5$ lines,
and hence they have the same combinatorial type.
(Recall that $\ell_{21}=\ell[56]$.)
On the other hand, we see that  $H(S_1)\sperp$ is an odd lattice and $H(S_2)\sperp$ is an even lattice.
Therefore they have different embedding topologies.
\begin{remark}
For the arrangement $S_1$, we have a line $\ell\sprime:=\ell_6$ satisfying $\intf{\ell\sprime, \ell\spprime}=0$ for any $\ell\spprime\in S_1$,
whereas there exists no such line $\ell\sprime$ for $S_2$.
\end{remark}
\begin{remark}
The fact that the set of all $5$-tuples of disjoint lines is decomposed into two orbits 
under the action of $W(E_6)$ was proved in~\cite[Proposition II-4]{Demazure1980}.
\end{remark}
\subsubsection{A pair of arrangements of $6$ lines}\label{subsubsec:6}
Let $o(T_1)$ and $o(T_2)$ be the $W(E_6)$-orbits whose minimal representatives are
\[
T_1=[1,2,3,4,5, 27]\quad{\rm and}\quad T_2=[1,2,3,4,21, 26], 
\]
respectively.
We have $|o(T_1)|=|o(T_2)|=432$.
Their combinatorial types are given by the dual graphs 
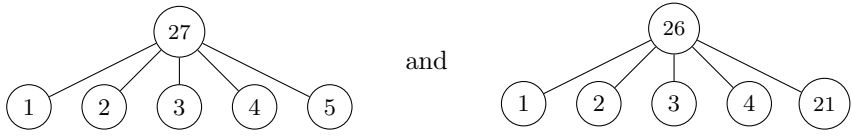
\[
\begin{tikzpicture}[x=1cm, y=1cm]
\node (L0)  [fill=white, draw, circle] at (0,1) {{\footnotesize $27$}};
\node (L1) [fill=white, draw, circle] at (-2,0) {{\small $1$}};
\node (L2) [fill=white, draw, circle] at (-1,0) {{\small $2$}};
\node (L3) [fill=white, draw, circle] at (0,0) {{\small $3$}};
\node (L4) [fill=white, draw, circle] at (1,0) {{\small $4$}};
\node (L5) [fill=white, draw, circle] at (2,0) {{\small $5$}};
\draw (L0)--(L1);
\draw (L0)--(L2);
\draw (L0)--(L3);
\draw (L0)--(L4);
\draw (L0)--(L5);
\end{tikzpicture}
\qquad\raise .8cm \hbox{\textrm{and}}\qquad 
\begin{tikzpicture}[x=1cm, y=1cm]
\node (L0)  [fill=white, draw, circle] at (0,1) {{\footnotesize $26$}};
\node (L1) [fill=white, draw, circle] at (-2,0) {{\small $1$}};
\node (L2) [fill=white, draw, circle] at (-1,0) {{\small $2$}};
\node (L3) [fill=white, draw, circle] at (0,0) {{\small $3$}};
\node (L4) [fill=white, draw, circle] at (1,0) {{\small $4$}};
\node (L5) [fill=white, draw, circle] at (2,0) {{\footnotesize $21$}};
\draw (L0)--(L1);
\draw (L0)--(L2);
\draw (L0)--(L3);
\draw (L0)--(L4);
\draw (L0)--(L5);
\end{tikzpicture},
\]
respectively, 
and hence they have the same combinatorial type.
(Recall that  $\ell_{21}=\ell[56]$, $\ell_{26}=\ell[\bar{5}]$ and $\ell_{27}=\ell[\bar{6}]$.)
We have
\[
H_1(X\setminus\Lambda(T_i))\cong H^3(X,\Lambda(T_i))\cong\Coker(H^2(X)\to \bigoplus_{\ell\in T_i}H^2(\ell))
\cong 
\begin{cases}
\ZZ/2\ZZ & \textrm{for $i=1$, }\\
0 & \textrm{for $i=2$.}\\
\end{cases}
\]
(Here we omit $\ZZ$ in the (co)homology groups.)
Therefore these two configurations  have different embedding topologies.
\bibliographystyle{plain}

\end{document}